\documentstyle[12pt]{article}
\titlepage

\newtheorem{thm}{Theorem}[section]
\newtheorem{lem}[thm]{Lemma}
\newtheorem{cor}[thm]{Corollary}

\begin{document}
\begin{titlepage}
\title{ On Adiabatic Limits and Rumin's Complex }
\author{Zhong Ge \thanks{
 Supported by the Ministry of Colleges and Universities
 of Ontario and the Natural Sciences and Engineering Research Council
of Canada }  \\
The Fields Institute
for Research in Mathematical Sciences \\
185 Columbia St. West, Waterloo, Ontario, N2L  5Z5 \\}
\date{ July 24, 1994}
\maketitle
\begin{abstract}
This paper shows that when the Riemannian metric on a contact manifold
is blown up along the direction orthogonal to the contact distribution,
the corresponding harmonic forms  rescaled and normalized
in the $L^2$-norms will converge to Rumin's harmonic forms.
This proves a conjecture in Gromov \cite{kn:gromov1}.
This result can also be reformulated in terms of spectral sequences, after
Forman, Mazzeo-Melrose. A key ingredient in the proof is the fact
 that the curvatures become unbounded in a controlled way.
\end{abstract}
\end{titlepage}

\section{Introduction}

Rumin  \cite{kn:rumin1}
constructed  a differential complex adapted to a contact distribution,
for which the Laplacians are sub-elliptic operators.
In this paper  we show how
to arrive at this complex via  adiabatic limits, using the ideas
of Mazzeo-Melrose \cite{kn:mazzeo}, and Witten \cite{kn:witten}.

Beginning with Witten's work on adiabatic limits \cite{kn:witten}, there is
a fair amount of work on  the asymptotic behaviors of geometric-topological
objects ( e.g.  harmonic forms, eta invariants, etc ) associated with
 a family of
Riemannian metrics on fiber bundles as the metrics
 become singular ( see, for example, Cheeger \cite{kn:cheeger}). In particular,
Mazzeo-Melrose \cite{kn:mazzeo}  studied those  of harmonic forms and related
them to spectral sequences ( see also Forman \cite{kn:forman}
). In all these work
an essential geometric assumption is that
 the curvatures of the metrics are uniformly  bounded.
In this
paper we consider a different   situation in which
a  Riemannian metric on a contact manifold
is blown up along the direction orthogonal to the contact distribution.
 It is known that,  despite that
 curvatures  become unbounded,
 the Riemannian metric  nevertheless
converges to a Carnot-Caratheodory metric,  and
  Gromov (\cite{kn:gromov1}, page 191-96 ) conjectured
 that  the harmonic forms will
converge to the corresponding objects associated with the
Carnot-Caratheodory metrics, i.e. the Rumin's harmonic forms.
In this paper we will show that this is indeed the case if the
harmonic forms are rescaled and normalized in the $L^2$-norms.
A key ingredient
in the proof is  the fact that the curvatures become unbounded in a controlled
way.

There is some interest to generalize  Rumin's theory to more general
Carnot-Caratheodory spaces ( see, for example,  Gromov \cite{kn:gromov1}).
 Some preliminary results
in this direction can be found in    \cite{kn:ge1}, \cite{kn:ge2}.
The results in this paper suggest
 that there is probably  a  different  approach, namely  that  through
the   study of   the adiabatic limits of
harmonic forms and the associated  `` spectral sequence '' $E^l_k$ ( cf. \S 2,
and Forman \cite{kn:forman}).
 This  is also related to
the characteristic cohomology ( cf.  Bryant-Griffiths \cite{kn:bryant} )
Vinogradov \cite{kn:vinogradov} ).

The results of this paper  have been announced in \cite{kn:ge4}.

\section{ Statement of Results}
 Let $M$ be a $(2m+1)$-dimensional compact Riemannian
manifold, $A$ a contact distribution. Let $B$ be
the orthogonal distribution to $A$,
so $TM=A \oplus B$.  Write the Riemannian metric as $g=g_A \oplus g_B$.
Consider
a family of metrics $g_\epsilon =g_A\oplus \epsilon ^{-2} g_B$. As $\epsilon
\to
0$, the metric space $(M, g_\epsilon )$
 converges in the sense of Gromov-Hausdorff to the Carnot-Caratheodory
metric space $(M, g_A)$ ( see, for example, Fukaya \cite{kn:fukaya},
Ge \cite{kn:ge3}, Gromov \cite{kn:gromov1} ),
 in which the distance between  two points
is the minimum of lengths of curves tangent to $A$ joining the two points.

Let $\Omega ^{p, q}= \Omega ^p(A) \wedge \Omega ^q(B)$. Decompose
$d$ into
$$ d=d^{2,-1} + d^{1, 0} + d^{0, 1}, \hspace{0.3in} d ^{a, b}: \Omega ^{p, q}
\to \Omega ^{p+a, q+b}.$$
Since $A$ is contact, $d^{2, -1}$ is not zero. This is the point
of departure of this paper  from Mazzeo-Melrose \cite{kn:mazzeo}.

Equip $\Omega ^{p, q}$ with the metric induced from $g_\epsilon
$ ( still to be denoted
by $g_\epsilon $ ).    Let $\Theta _\epsilon $ be the isometry ( the rescaling
map )
$$ \Theta _\epsilon : ( \Omega (M), g_\epsilon ) \to ( \Omega (M), g_1). $$
Define the normalized differential $d_\epsilon =\Theta _\epsilon \circ d  \circ
(\Theta _\epsilon)^{-1} $, then
$$ d_\epsilon = {1\over \epsilon } d^{2, -1} + d^{1, 0} + \epsilon
 \;  d^{0, 1}.$$
We will use  `` $*$ '' to   denote the adjoint with respect to $g_1$.

Now Rumin's complex ( see Rumin \cite{kn:rumin1} )  can be written as
$${\cal R}^k := \Omega ^{k, 0} /Im(d^{2, -1}),\; \hspace{0.1in}  k\leq m;
\hspace{0.2in} {\cal R}^k := \Omega ^{k-1, 1} \bigcap Ker(d^{2, -1}),
  k \geq  m+1, $$
with the induced differential
\begin{eqnarray*}
d_\xi & = &\pi \; d ^{1, 0}:\hspace{0.2in} {\cal R}^k \to {\cal R}^{k+1},\;
\hspace{1.4in}
 k \ne m; \\
 d_{{\cal R}} & = &\pi (d^{0, 1}-d^{1, 0}
(d ^{2, -1} ) ^{-1} d^{1, 0}): {\cal R}^m \to {\cal R}^{m+1},
\end{eqnarray*}
where $\pi $ is the orthogonal projection $\Omega ^k\to {\cal R}^k$.
This is a sub-elliptic complex.
\begin{thm} Assume  ($M, g$ ) is 
 Heisenberg  ( cf. \S 3).
 Suppose $\omega _\epsilon, \Vert \omega _\epsilon \Vert _{L^2} =1$,
is a $d_\epsilon$-harmonic form,
$$ d_\epsilon \omega _\epsilon =d^*_\epsilon \omega _\epsilon=0,  $$
(i.e. $\Theta ^{-1}_\epsilon \omega _\epsilon $ is a harmonic form for $(M,
g_\epsilon)$ ).  Then as
$\epsilon \to 0$,  after
passing to a subsequence,
$$ \omega _\epsilon \to \omega _0 \ne 0
 \hspace{0.3in}\mbox{ strongly in $L^2$,}
$$
 and $\omega _0$ is a Rumin's harmonic form.
\end{thm}

This result can also be reformulated in terms of spectral sequence,
after Mazzeo-Melrose \cite{kn:mazzeo}, Forman \cite{kn:forman}.

Fix a number $l$, one says a family of $k$-forms $\omega _\epsilon $ is of
class
$O(\epsilon ^l)$, i.e. $ f_\delta = O(\epsilon ^l) $,
if $ \epsilon ^{-l}    \Vert \omega _\epsilon \Vert
_{H^1}$ is uniformly bounded, where $H^1$ denotes the ordinary Sobolev space.
Define
$$ E^k_l := \{ \;  \omega _\epsilon \in \Omega ^k,
d_\epsilon \; \omega _\epsilon =
O(\epsilon ^{l-1}), \; \; d_\epsilon ^*  \; \omega _\epsilon =O(\epsilon
^{l-1}),
\; \; \; \Vert \omega _\epsilon \Vert _{L^2}=1 \},  $$
and set
$$ \bar E^k_l := \mbox{ linear span  }
 \{ \mbox{ the weakly limits of
$ \omega _\epsilon $ in $L^2 $ as $\epsilon \to 0 $ },
\; \omega _\epsilon \in E^k_l \; \} \;  \bigcap C^\infty (\Omega ^k).$$
Obviously, for each $k$
$$ \cdots \cdots  \subset \bar E^k_2 \subset  \bar E^k_1 \subset
 \bar E^k_0.$$
The following result
says that most terms in the spectral sequence will degenerate at
$\bar E_2$ except those of degree either $m$ or $m+1$, which degenerate at
$\bar E_3$.
This may explain
 why $d_{{\cal R}}$ is a second order operator.
\begin{thm}  Suppose ($M, g$ ) is 
Heisenberg ( cf. \S 3).

 (1) The terms in $\bar E^1$ are
$$   \bar E^k_1= {\cal R}^k.$$
(2) The terms in $\bar E^2 $ are
\begin{eqnarray*}
 &&  \bar E^k_1=\{ \omega \in {\cal R}^k,  d _\xi \omega =d_\xi ^* \omega =0
\},
   \hspace{0.4in} k \ne m, m+1; \\
 &&  \bar E^k_2 = \{ \omega \in {\cal R}^k,  ( d^{1, 0} )^* \omega =0 \},
\hspace{0.6in} k =m; \\
 &&  \bar E^k_2 = \{ \omega \in {\cal R}^k,  d^{1, 0} \omega =0 \},
\hspace{0.8in} k =m+1.
\end{eqnarray*}
(3) The terms in $\bar E^3$ are
\begin{eqnarray*}
 &&  \bar E^k_3= \bar E^k_2,
   \hspace{2.65in} k \ne m, m+1; \\
  && \bar E^k_3 =\bar E^k_4=\cdots =
 \{ \omega \in {\cal R}^k, \; \; d _{{\cal R}} \omega =d _\xi ^* \omega =0 \}
, \hspace{0.1in} k=m; \\
 &&  \bar E^k_3 =\bar E^k_4=\cdots = \{ \omega \in {\cal R}^k, \; \; d _\xi
\omega =d _{{\cal R}}
 ^* \omega =0 \}, \hspace{0.1in} k=m+1.
\end{eqnarray*}
\end{thm}

We shall use the following notations: $\Vert \cdot \Vert _{H^1_c}$
denotes the  following weighted Sobolev's norm
$$ \Vert \omega \Vert ^2_{H^1_c}=\int _M  \sum
(D_{e_i} \omega , D_{e_i} \omega )
dv, $$
where $e_i$ is an orthonormal basis for $A$, and
$\Vert \omega \Vert _{H^2_c}$ similarly.

\section{  Geometry of  Heisenberg Manifolds }

Let $v$  be
a ( local ) unit tangent vector field
 spanning $B$, and $\xi$ the contact 1-form which
satisfies $\xi (v)=1$.  If
 the metric $g_A$ can be written as $ g_A(a, b )= d\xi (a, J b), a, b \in A$,
where $J$ is an endmorphism of $A$ satisfying $J^2=-Id $,
then we say ( $M, g$ )
 is { \bf 
 Heisenberg. } Note that even though $v$
is in general only locally defined, this notion  is well defined.
Throughout the
rest of this paper we assume that ( $M, g$ ) is  Heisenberg.

We will use the following properties  of
Heisenberg manifold.
\begin{lem}
There is an orthonormal basis $e_1, \cdots, e_m, e_{m+1}:= J e_1,
\cdots, e_{2m}:=J e_{m}$ for $A$ such that
\begin{eqnarray*}
&& [ e_i, e_j ]  =  0 \hspace{0.5in}   \mbox{mod} \,\, (A), \nonumber \\
 &&  [ e_i, e_{m+j} ] =  \delta _{i j} \; v \hspace{0.1in} \mbox{mod} \, \,
(A),
\hspace{0.2in} 1\leq i, j \leq m.
\end{eqnarray*}
\label{le:1heisen}
\end{lem}
\noindent {\bf Proof. }  
This follows from the condition
$ g_A=d\xi (\cdot, J \cdot ), $  $J^2=-Id. $
\vspace{0.1in}

The following property in fact holds for any contact manifold.
\begin{lem}
Suppose  $\xi $ is a contact $1$-form and $x_0$ a fixed point on $M$,
 $v$  the vector field such that $i(v)\xi=1, i(v) d\xi=0$.
 There are vector fields
$u_i$, $i=1, \cdots, 2m+1$, $u_{2m+1}=v$, such that
\begin{enumerate}
\item $u_i $ are linearly independent at $x_0$;

\item $u_i$  vanishes  outsider a
small neighborhood for $i=1,\cdots, 2m$;

\item ${\cal L}_{u_i} \xi ={\cal L}_{u_i} d\xi=0$, and ${\cal L}_{u_i}$
 preserves $\Omega ^{k, 0}$, $i=1, \cdots, 2m+1$.
( Here ${ \cal L }_u$ is the Lie derivative in the direction $u$.)
\end{enumerate}
 \label{le:contact}
\end{lem}
\noindent {\bf Proof. } First note that ${\cal L}_v \xi =0$ follows from
${\cal L}_v=i(v) d +d\; i(v)$.

To choose $u_1, \cdots, u_{2m}$, one takes a local coordinates $(x, z)\in {\bf
 R^{2n}}\times {\bf R}$
near $x_0$ such that
$ \xi = dz- \rho $ where $\rho $ is a $1$-form on ${\bf R^{2n}}\{x\}$ and
$v=\partial / \partial z$.
Choose $2m$  functions $f_1, \cdots, f_{2m}$ on ${\bf R^{2m}}\{x\}$
  with  linearly independent
$df_1, \cdots, df_{2m}$ at $x_0$ such that
 $f_i$  vanishes outsider a neighborhood.
Let $H_{f_i}$ denote the Hamiltonian vector
field of  $f_i$  with respect to
$d\rho $.
Define $u_i (x, z)=H_{f_i} (x)+  a_i \partial / \partial z$,  where
$a_i$ is determined from the equation $i(u_i) \xi =-f_i$,
$i=1, \cdots, 2m$. One easily verifies  that $u_i$ thus defined satisfies
all the  requirements.
\vspace{0.1in}

\section{ $ A\;priori$ Estimates }

To prove Theorem 2.1  and Theorem 2.2,
in this section we will derive some $ a\;priori$ estimates for
 the $H^1_{c}$-norm
of $\omega$ in terms of $ ( \Delta _{d _\epsilon } \omega, \omega )$ if
 $k\ne
m, m+1$,   and for  the $H^2_c$-norm of $\omega $ if $k=m, m+1$.  As the case
of $k>m$ is similar to  that of $k\leq m$, we will only consider
the  case $ k\leq m$.

We will use the following notations: If $L$  is an operator,
$$ \Delta _{L} := L^*\;  L +L\;  L ^*. $$
The  letter $C$ denotes a generic positive number, $M$
a generic constant.

\subsection{ The case of $k$-forms ( $k \ne m, m+1$ ).}

We have  the following  $a\;priori$ estimates
\begin{thm}
For any $\omega =\alpha+ \beta,  \alpha \in \Omega ^{k, 0},
\beta \in \Omega ^{k-1, 1}$, we have
\begin{eqnarray*}
( \Delta _{d_\epsilon } \omega, \alpha  ) \geq  {1 \over \epsilon
^2} \Vert (d^{2, -1})^*
 \alpha \Vert ^2_{L^2} +
 {m-k \over m } C
\Vert \alpha \Vert _{H^1_c} ^2 + \epsilon ^2( D_v \alpha , D_v \alpha  )
& - & M (\omega, \omega ),  \\
&& k\leq m-1;
\end{eqnarray*}
and
\begin{eqnarray*}
( \Delta _{d_\epsilon } \omega, \beta ) \geq {1 \over \epsilon ^2}
\Vert (d^{2, -1})
 \beta \Vert ^2_{L^2} + {m-k+ 1 \over m } C
\Vert \beta \Vert _{H^1_c} ^2 + \epsilon ^2 ( D_v \beta , D_v \beta )
& - & M (\omega, \omega ), \\
& & k\leq m.
\end{eqnarray*}
 \label{le:apr}
\end{thm}
To  prove this theorem, we need a few technical results.
\begin{lem} The following operator
$$ Q: = (d ^{0, 1})^* d^{1, 0} + (d^{1, 0}) ^* d^{ 0, 1}
+ d^{0, 1} ( d^{1, 0})^* +d^{1, 0} (d^{0, 1})^*$$
is a first-order linear differential operator.
\end{lem}
\noindent {\bf Proof. }  We only need to prove that
$$ (d ^{0, 1})^* d^{1, 0} + (d^{1, 0}) (d^{ 0, 1} )^* $$
is a first order operator.

 If $e_i$ is an orthonormal basis for $A$, $v$ for $B$, then
 \begin{eqnarray*}
& &  d^{1, 0}  =\sum e^i \wedge D_{e_i} +\mbox{ 0-order operator }, \\
 & & (d^{0, 1})^*   = i(v) D_v + \mbox{ 0-order operator }.
\end{eqnarray*}
So
\begin{eqnarray*}
& & ( d ^{0, 1} )^* d^{1, 0} + d^{1, 0}  ( d ^{0, 1} )^* =\\
&& =\sum e^i \wedge D_{e_i}  i(v) D_v + i(v) D_v e^i \wedge D_{e_i}  +
\mbox{ 1st order operator } \\
&& = \sum e^i \wedge i(v)
 D_{e_i}   D_v + i(v)  e^i \wedge D_{e_i}D_v   +
\mbox{ 1st order operator } \\
& &= \mbox{ 1st order operator. }
\end{eqnarray*}
Here we have used the fact that $e^i \wedge i(v)+ i(v) e^i \wedge  =0$.
\vspace{0.1in}
\begin{lem} If $\alpha, \beta $ are as in Theorem \ref{le:apr}, then
\begin{eqnarray*}
&& (\Delta _{(\epsilon ^{-1} d^{2, -1} +d ^{1, 0})}\; \alpha,
\alpha )\geq {1 \over \epsilon ^2}
 \Vert (d^{2, -1})^* \alpha \Vert ^2 _{L^2} + C {m-k \over m} \Vert \alpha
\Vert _{H^1_c}^2 -M( \alpha, \alpha ), \; k\leq m-1; \\
& & (\Delta _{ (\epsilon ^{-1} d^{2, -1} +d ^{1, 0} )}\; \beta,
\beta )\geq
{1 \over \epsilon ^2}
 \Vert d^{2, -1}\beta  \Vert ^2 _{L^2}+
C {m-k+ 1 \over m} \Vert \beta \Vert _{H^1_c}^2 -M( \beta, \beta ),
 \; k\leq m. \end{eqnarray*}
 	\label{le:le}
\end{lem}
\noindent {\bf Proof. } We shall only prove the first inequality,
as the second one can be proved similarly.

 By counting the types of the differential forms, one
has
$$ (\Delta _{ (\epsilon ^{-1} d^{2, -1} +d ^{1, 0} )}\hspace{0.1in} \alpha,
\alpha )\geq {1 \over \epsilon ^2}
 \Vert (d^{2,1})^* \alpha \Vert ^2 + ( \Delta _{d ^{1, 0}} \;
 \alpha ,  \alpha ). $$
In terms of a local orthonormal
basis $e_i$ for $A$, one can write
\begin{eqnarray*}
 d^{1, 0} & = & \sum _{i=1}^{2m} e ^i \wedge D_{e_i}+
\mbox{ 0-order operator }, \\
(d^{1, 0})^* & = &  \sum _{i=1}^{2m} i(e_i) \wedge D_{e_i} + \mbox{ 0-order
operator }.
\end{eqnarray*}
So
\begin{eqnarray*} \Delta _{d^{1, 0}}  &= &
\sum e^i \wedge i(e_j) D_{e_i} D_{e_j} +i(e_j) e^i \wedge D_{e_j} D_{e_i}
+\mbox { 1-st order operator in $e_i$ } \\
& =& -\sum  D_{e_i} D_{e_i} -e^i \wedge i(e_j)  ( D_{e_j} D_{e_i}-
D_{e_i} D_{e_j} )
+\mbox { 1-st order operators in $e_i$ }.
\end{eqnarray*}
Now by Lemma \ref{le:1heisen}, if $i>j$,
$$ D_{e_i} D_{e_j}-D_{e_j} D_{e_i} ={1\over m}  \delta _{i-m, j}
\sum _{i=1}^m ( D_{e_{i+m}} D_{e_i} -D_{e_i} D_{e_{i+m}} )
+\mbox{ 1st order operator in $D_{e_i}$ },$$
while if $i<j$,
$$ D_{e_i} D_{e_j}-D_{e_j} D_{e_i} = -{1\over m} \delta _{i+m, j}
\sum _{i=1}^m ( D_{e_{i+m}} D_{e_i} -D_{e_i} D_{e_{i+m}} )
+\mbox{ 1st order operator in $D_{e_i}$ }.$$
So after an integration by parts,
\begin{eqnarray}
&& \int_M  ( e^i \wedge i(e_j) ( D_{e_i} D_{e_j} \alpha -
D_{e_j} D_{e_i} \alpha, \alpha ) \nonumber \\
&& ={2 \over m}
 \int _M (i( e_i)  D_{e_{j+m}} \alpha,  i(e_{i+m}) D_{e_j}  \alpha )
-(i( e_i)  D_{e_j} \alpha,  i(e_{i+m}) D_{e_{j+m}}  \alpha ) \nonumber \\
&& + \mbox{ terms
of the form $ ( D_{ e_i} \alpha , \alpha ) $ } \nonumber \\
&& \leq  {k\over m} \Vert \alpha \Vert _{H^1_c}^2
+\mbox{terms
of the form $ (  D_{e_i} \alpha, \alpha ) $. }  \label{eq:key}
\end{eqnarray}
This proves the lemma.
\vspace{0.1in}

\noindent {\bf Proof of Theorem \ref{le:apr}.}
 We shall only prove the first inequality,
as the second one can be proved similarly.

 By a direct computation, one has
\begin{eqnarray}
  (\Delta \omega, \alpha ) & =& (D_{ (\epsilon ^{-1} d^{2, -1} +d ^{1, 0} )}
\hspace{0.1in} \omega, \alpha ) + \label{eq:1st} \\
& & + ( ( d^{2, -1})^* \alpha, (d^{0, 1})^* \omega )
+ ( ( d^{2, -1})^* \omega , (d^{0, 1})^* \alpha ) \label{eq:2nd1}\\
& & + ( ( d^{2, -1})^* \alpha, (d^{0, 1})^* \omega )
+ ( ( d^{2, -1})^* \alpha, (d^{0, 1})^* \omega ) \label{eq:2nd2}\\
& & + \epsilon
 ( Q \omega, \alpha ) +\epsilon ^2 ( \Delta _{d^{0, 1}}
 \hspace{0.1in}
\omega, \alpha ).  \label{eq:3rd}
\end{eqnarray}
The first term (\ref{eq:1st})  was considered in Lemma \ref{le:le}.
By counting the types, the terms (\ref{eq:2nd1}), (\ref{eq:2nd2})
are zero. We only need to treat
the remaining term (\ref{eq:3rd}).

Note that
$$ (\Delta _{d^{0, 1}} \hspace{0.1in} \omega, \alpha ) = (D_v \alpha,
D_v\alpha )+ \int ( \mbox{ terms of the form $( D_v \omega, \alpha  )$})$$
By the Schwartz inequality and the fact that $
Q$ ( and $Q^*$ ) is a first order operator,
$$
 \epsilon
 (Q\omega, \alpha )= \epsilon
 (\omega, Q^* \alpha )\;  > \;  - {\epsilon ^2  \over 2}
 (\Vert \alpha \Vert_{H^1_c
}^2 + \Vert D_v \alpha \Vert _{L^2}^2 ) - M \Vert \omega \Vert _{L^2}^2.
$$
Substituting these inequalities into eq. (\ref{eq:1st})-(\ref{eq:3rd}), we
prove the theorem.
\vspace{0.2in}

\subsection{ The case of $m, (m+1)$-forms. }
If $k=m,$  the estimate for the derivatives  of
$\alpha $ in Theorem \ref{le:apr} breaks down.
So we need  a different method to do the  estimate.

Suppose $\omega =\alpha +\beta, \alpha \in \Omega ^{m, 0},
\beta \in \Omega^{m-1, 1}$,   satisfies
\begin{eqnarray}  & & d_\epsilon \omega =\xi _1, \label{eq:xi} \\
& & d^*_\epsilon \omega = \xi _2. \label{eq:label}
\end{eqnarray}
Let $\gamma _1 = D_v \alpha, \gamma _2=D_v \beta $, $\gamma =\gamma _1
+\gamma _2$.
We first estimate the second order derivatives of $\beta. $

Note that $D_v$ commutes with $d^{2, -1}$ and
with  $(d^{2, -1})^*$  modulo  zero-order operators.
Thus,
 taking the derivative $D_v$ of the eqs. (\ref{eq:xi}),  (\ref{eq:label}),
one obtains
\begin{eqnarray*}
 d_{\epsilon } \;  \gamma & = & D_v \xi _1
+  (\mbox{  0-order operator in powers of   $\epsilon $   } )\;
 \omega \\
& & +  {1\over  \epsilon } (\mbox{  0-order operator    } ) \; \omega,
\\
  d_{\epsilon } ^*  \;  \gamma & = &
D_v \xi _2
+  (\mbox{  0-order operator in powers of   $\epsilon $   } ) \; \omega \\
& & + {1\over  \epsilon } (\mbox{  0-order operator    } )
\; \omega.
\end{eqnarray*}
Integrating
$$\int _M (\Delta _{d_{\epsilon }} \;  \gamma, \gamma _2 ).
$$
as in the proof of Theorem \ref{le:apr},  using
 the facts that $\omega, \gamma _2$ are uniformly bounded in the
$L^2$-norm ( Theorem \ref{le:apr} ), one obtains
\begin{lem} Suppose $k=m$.
If $\omega =\alpha +\beta $ satisfies eqs. (\ref{eq:xi}), (\ref{eq:label}),
then
$$ \Vert D_v \beta \Vert _{H^1_c} ^2\leq C (
{ 1\over \epsilon ^2}  \Vert \omega \Vert _{L^2} ^2
+ \Vert \beta  \Vert ^2 _{L^2} + \Vert D_v \xi \Vert _{L^2}^2 ).
$$ \label{le:initial}
\end{lem}

Now we estimate $\Vert \alpha \Vert _{H^2_c } $.  For this  purpose
we need to  decompose $\alpha =\alpha _1 +L\alpha _2 $,
  since $\Omega ^{m, 0}=
{\cal R}^m \oplus L(\Omega ^{m-2, 0})$,
 where $\alpha_1 \in {\cal R}^m$,  $L: \Omega ^{m-2, 0}\to
\Omega ^m $ is defined by $L \alpha _2=\alpha _2 \wedge d\xi$, and
$\alpha _1 $, $L \alpha_2 $ are orthogonal.  We will  estimate the
derivatives of $\alpha _1 $ and $L\alpha _2$ separately.
We first estimate the first order derivatives of $L \alpha_2$.
\begin{lem}
$$\Vert L \alpha _2  \Vert _{H^1_c} ^2 +
 { 1\over \epsilon ^2}  \Vert (d^{2, -1})^* L\alpha _2 \Vert _{L^2}^2 \leq C (
\Vert \xi  \Vert _{L^2}^2 +\Vert \omega \Vert
_{L^2} ^2 ).$$
  \label{le:we}
\end{lem}
\noindent {\bf Proof. } This is proved by integrating
$$ \int ( L\alpha _2, \Delta _{d_\epsilon  }\; \omega )$$
as in the proof of Theorem \ref{le:apr}. The key point is that
$\Delta _{d^{1, 0}}$ is sub-elliptic on $L(\Omega ^{m-2, 0})$. In fact,
one has the following estimate which improves over that in (\ref{eq:key})
\begin{eqnarray*}
&& \int_M  ( e^i \wedge i(e_j) ( D_{e_i} D_{e_j}  -
D_{e_j} D_{e_i}) L\alpha _2, L\alpha _2) \nonumber \\
&& ={2 \over m}
 \int _M (i( e_i)  D_{e_{j+m}} L\alpha_2,  i(e_{i+m}) D_{e_j}  L\alpha_2 )
-(i( e_i)  D_{e_j} L \alpha _2, i(e_{i+m}) D_{e_{j+m}} L \alpha_2 ) \nonumber
\\
&& + \mbox{ terms
of the form $ (  D_{e_i} L \alpha _2 , L \alpha_2 ) $ } \nonumber \\
&& \leq  {m-1\over m} \Vert L \alpha _2 \Vert _{H^1_c}^2
+\mbox{terms
of the form $ ( D_{ e_i} (L\alpha _2 ), L \alpha _2 ) $ },
\end{eqnarray*}
and hence  one has
$$ (\Delta _{d_\epsilon}\;  \omega, L\alpha _2) \geq
 {1 \over \epsilon ^2}
\Vert (d^{2, -1})^*
L\alpha _2 \Vert ^2_{L^2} + { 1 \over m } C
\Vert L\alpha_2 \Vert _{H^1_c} ^2
 - M (\omega, \omega ). $$
\vspace{0.2in}

We now  estimate the second order derivatives of $L\alpha _2$. Let
$u=u_i$ be as in Lemma \ref{le:contact}, $i=1, 2, \cdots, 2m+1$. Take the Lie
derivative of
eqs. (\ref{eq:xi}),  (\ref{eq:label} ) with respect to $u$, one has
\begin{eqnarray*}
 d_{\epsilon } \;  \varphi & = & {\cal  L}_{u}\;  \xi _1
+  (\mbox{  0-order operator in powers of   $\epsilon $   } )\;
 \omega, \\
  d_{\epsilon } ^*  \;  \varphi & = &
{\cal L}_u\;  \xi _2
+  (\mbox{  0-order operator in powers of   $\epsilon $   } ) \; \omega,
\end{eqnarray*}
where $\varphi := {\cal  L}_u \omega$. Note that $
{\cal L}_{u_i}$ preserves the decomposition
$\Omega ^{m, 0}={\cal R}^m \oplus L(\Omega ^{m-2, 0})$,
 hence,  applying the same arguments as in the proof of
 Lemma \ref{le:we} to the
above equations, one obtains
\begin{lem}
$$ \Vert {\cal L}_{u_i} L \alpha _2  \Vert _{H^1_c}^2 \leq C (
\Vert {\cal  L}_{u_i} \xi  \Vert _{L^2}^2   + \Vert \xi  \Vert _{L^2}^2  +
\Vert
\omega  \Vert
_{L^2} ^2 ), \hspace{0.2in} i=1, 2, \cdots, 2m+1.
$$ \label{le:fen}
\end{lem}

Then we have the following estimate on the second order derivatives of
$L\alpha _2$.
\begin{cor}
$$ \sum _{i=1}^{2m}\Vert D_{e_i}  L\alpha _2  \Vert _{H^1} ^2\leq C (
\Vert  \xi  \Vert _{H^1}^2   +  \Vert  \xi  \Vert
_{L^2} ^2 +  \Vert  \omega  \Vert
_{L^2} ^2 ).
$$
\label{le:cor5}
\end{cor}
\noindent {\bf Proof. } This follows
from Lemma \ref{le:fen}. In fact,
by Lmmma \ref{le:fen},  every $x_0\in M$ has a neighborhood $U$ such that
$$ \sum _{i=1}^{2m}\Vert D_{e_i}  L\alpha _2  \Vert _{H^1(U)} ^2\leq C (
\Vert  \xi  \Vert _{H^1}^2   +  \Vert  \xi  \Vert
_{L^2} ^2 +  \Vert  \omega  \Vert
_{L^2} ^2 +\Vert L \alpha _2\Vert _{H^1_c}).
$$
Then the corollary follows from  a partition of unity and Lemma \ref{le:we}.
\vspace{0.2in}

Now we estimate the derivatives of $\alpha _1$.

Eliminating $\beta $ from eq.  (\ref{eq:xi}) by  the fact that
$ d^{2, -1}: \Omega ^{m-2, 1} \to \Omega ^{m+1, 0}$ is an isomorphism, one
has
\begin{equation}
 d_{{\cal R}} \alpha_1 = \pi (  ( d^{2, -1})^{-1} \; \xi^1_1 +  { 1\over
\epsilon }
\xi _1^2 )-d_{{\cal R}} \; L \alpha _2,
\label{eq:pr1}
\end{equation}
where $\xi _1 = \xi ^1_1 + \xi ^2 _1, \xi ^i \in \Omega ^{m-i-1, i-1}, i=1, 2$.
{}From the $(1, 0)$-component of eq. (\ref{eq:label}) one  obtains
\begin{equation}
 d ^{1, 0} (d ^{1, 0})^* \alpha _1
=   d^{1, 0} \; ( - \epsilon \; ( d^{0, 1}) ^* \; \beta  + \xi _2^1 )-
 d ^{1, 0} (d ^{1, 0})^* L \alpha _2,
\label{eq:pr2}
\end{equation}
where $\xi _2$ is decomposed into $\xi ^1_2 + \xi _2^2$ as for
$ \xi _1$.
By  Rumin \cite{kn:rumin2},
 $ d_{{\cal R}} ^* d_{{\cal R}} + (d^{ 1, 0} ( d^{1, 0}) ^*)^2
$ is hypoelliptic on ${\cal R}^m$.
(  However, note that $ d_{{\cal R}} ^* d_{{\cal R}} + (d^{ 1, 0} ( d^{1, 0})
^*)^2
$ is not hypoelliptic  on $L(\Omega ^{m-2, 0})$, which is the
reason why we decompose $\alpha $. )   Hence,  from the eqs.
 (\ref{eq:pr1})-(\ref{eq:pr2}), plus  the following estimates
\begin{eqnarray*}
&& \Vert d_{{\cal R}} L\alpha _2  \Vert _{L^2}\leq C  (
\sum _{i=1}^{2m} \Vert D_{e_i} L\alpha _2 \Vert _{H^1} +  \Vert
 L\alpha _2 \Vert _{L^2}), \\
&& \Vert d^{1, 0} (d^{1, 0})^* L\alpha _2  \Vert _{L^2}\leq C (
\sum _{i=1}^{2m} \Vert D_{e_i} L\alpha _2 \Vert _{H^1}+\Vert
 L\alpha _2 \Vert _{L^2}),
\end{eqnarray*}
which can be controlled by using Corollary
\ref{le:cor5},   and
$$ \Vert  d^{1, 0} ( d^{0, 1}) ^* \; \beta \Vert _{L^2}
\leq C ( \Vert D_v \beta \Vert ^2 _{H^1_c}  + \Vert \beta \Vert ^2 _{H^1_c}
+\Vert \beta \Vert ^2_{L^2}),
$$
which can be controlled by using  Lemma \ref{le:initial}, one obtains
\begin{thm} If $\omega =\alpha + \beta $ satisfies eqs. (\ref{eq:xi}),
(\ref{eq:label}), then
$$ \Vert \alpha \Vert _{H^2_c}^2 \leq
 C ( \Vert \omega \Vert _{L^2}^2
 + \Vert \xi \Vert ^2 _{H^1}+ {1 \over \epsilon ^2}
\Vert \xi \Vert ^2 _{L^2} ).
$$ \label{le:apr2}
\end{thm}
\section{ Proof of Theorem 2.1. }

Much of the proof depends on the properties of $d^{2, -1}$, which we list
now. These  properties follow from a straight forward computation.
\begin{lem} (1) $ d^{2, -1}: \Omega ^{k-2, 1} \to \Omega ^{k, 0}$ is an
injection for $k\leq m-1$, an
isomorphism for $k=m$.

(2) $ (d ^{2, -1} )^* : \Omega ^{k, 0} \to \Omega ^{k-2, 1}$ is an injection
for $k\geq m+2 $, an isomorphism for $k=m+1$.
\label{le:d21}
\end{lem}
Let $\omega =\alpha + \beta $, $\alpha \in \Omega ^{k, 0}, \beta \in
\Omega ^{k-1, 1}$ be as in Theorem 2.1. Then $d_\epsilon \omega =( d_\epsilon )
^*\omega =0 $ is equivalent to
\begin{eqnarray}
&&  { 1\over \epsilon } d^{2, -1} \beta + d^{1, 0} \alpha =0, \label{eq:eq1}
 \\
&& d^{1, 0} \beta + \epsilon  \; d^ {0, 1} \alpha =0, \label{eq:eq2} \\
&& { 1\over \epsilon } (d^{2, -1})^* \alpha + (d^{1, 0})^* \beta =0,
\label{eq:eq3}
 \\
&& ( d^{1, 0} )^*\alpha + \epsilon ( d^{0, 1})^* \beta =0. \label{eq:eq4}
\end{eqnarray}

\begin{lem} Suppose ${\cal Q}$ is a first-order differential operator.
 If  $\Vert \omega _\epsilon  \Vert _{L^2}$,  $\Vert
\epsilon \; {\cal Q} \omega_\epsilon\Vert _{L^2}$ are uniformly bounded,
then
$$ \epsilon\; { \cal Q}\omega _\epsilon \to 0, \hspace{0.2in}
\mbox{ weakly in } L^2. $$ \label{le:cor}
\end{lem}
\noindent {\bf Proof. } We may  choose a sub-sequence such that
$$ \epsilon \;  {\cal  Q } \omega_\epsilon  \to \gamma, \hspace{0.2in}  \mbox{
weakly in } L^2. $$
We need only to prove $\gamma =0$. Assume without loss of generality
that $ \omega \to a $ weakly in  $L^2$.  Now choose a smooth $k$-form $\mu$,
then
$$ ( \gamma, \mu )_{L^2}
 = lim_{\epsilon \to 0}\;  ( \epsilon \; {\cal Q }\omega _\epsilon, \mu )_{L^2}
= lim_{\epsilon \to 0 } \epsilon
\;   ( \omega _\epsilon , ({\cal Q})^* \mu )_{L^2}=0,$$
so $\gamma =0$.
\vspace{0.2in}

\noindent {\bf Proof of Theorem 2.1.}
We will only prove the theorem for the case
 $k\leq m$, as the case $k>m$ is similar.
We divide the proof into two cases: one for $k\leq m-1$, the other $k=m$.

(1) $k\leq m-1$.  By  Theorem \ref{le:apr} we observe that both $\alpha $
and $\beta $ are
uniformly  bounded in $H^1_{c}$, and $\epsilon ^{-1} d^{2, -1} \beta $ is
uniformly bounded
in $L^2$.  By Lemma 5.1, this implies that $\alpha $ converges
to $\alpha _0$, after passing to a subsequence if necessary, and
$\beta $  to $0$
strongly in $L^2$.

By Theorem \ref{le:apr}, $\epsilon
\; \Vert d^{0, 1} \beta \Vert _{L^2}$ is bounded.
Then by Lemma \ref{le:cor}, $\epsilon \; d^{0, 1} \beta \to 0$ weakly
in $L^2$. Then,
from eq. (\ref{eq:eq2}),  it follows that $\alpha _0$
satisfies
$$ d^{1, 0} \alpha _0 =0$$
in the weak sense that for any $\mu \in H^1_c$,
$( \alpha _0,  (d^{1, 0})^*\; \mu )_{L^2}=0$.

Similarly, from eq. (\ref{eq:eq3}) and Lemma 5.2 we have
$$ (d^{1,0})^* \alpha _0 =0, \hspace{0.2in} (d^{2, -1}) ^* \alpha _0=0$$
in the  weak sense. Now the theory of sub-elliptic  operators implies
that $\alpha _0$ is smooth and satisfies the Rumin's Laplacian ( cf.
 Helffer-Nourrigat \cite{kn:helffer} ).

To conclude the proof,
 we note that $\alpha _0 \ne 0$, as $\omega $ converges to $\alpha _0$ strongly
in $L^2$, and $\Vert \omega \Vert _{L^2}=1$.
  This proves the theorem for $k< m$.

(2) If $k=m$, then it follows from Theorem \ref{le:apr},
 Theorem \ref{le:apr2},  that
$\Vert \alpha \Vert _{H^2_{c}} $, $\Vert \beta\Vert _{
H^1_c}$ are uniformly bounded.
 Moreover, as in the case $k\leq m-1$, $\beta \to 0$ weakly
in $H^1_c$.
We may choose a subsequence of $\alpha $ such that $\alpha \to
\alpha _0$ weakly in $H^2 _c$.

It follows from
 eqs. (\ref{eq:eq1}), (\ref{eq:eq2}),  that
$$ d_{\cal R} \alpha =0.$$
Hence $\alpha _0$ also satisfies the above equation in the weak sense that
for any $\mu \in H^2_c$, $(\alpha _0, (d_{{\cal R}})^*\mu )_{L^2}=0.$

Now that $\epsilon \;  ( d^{1, 0})^*\beta,
\epsilon \;  ( d^{0, 1})^*\beta
$ are uniformly bounded in $L^2$,
by   Lemma \ref{le:cor},
   $ \epsilon \;  ( d^{1, 0})^*\beta \to 0$,
$\epsilon \;  ( d^{0, 1})^*\beta  \to 0 $ weakly in $L^2$.
Then it follows from eqs. (\ref{eq:eq3}), (\ref{eq:eq4}) that
$$ ( d^{1, 0} ) ^* \alpha _0 =0,
 \hspace{0.2in} (d^{2, -1}) ^* \alpha _0=0
$$
in the weak sense. Thus $\alpha _0$ is a Rumin's harmonic form.

At last  note that since $\omega \to \alpha _0$ strongly in $L^2, $
 $ \alpha _0 \ne 0$. This proves the theorem.
\vspace{0.2in}

\section{ Proof of Theorem 2.2. }

Ae before,  we only
consider the case $k\leq m$,
as the case  $k>m $ is similar.

First note that
 the equations $ d_\epsilon \omega =\xi_1, d _\epsilon \omega =\xi _2$,
are equivalent to
\begin{eqnarray}
  && {1 \over \epsilon } d^{2, -1} \beta + d^{1, 0} \alpha =\xi _1^1,
\label{eq:eq5} \\
  && d^{1, 0}\beta + \epsilon \; d^{0, 1 } \alpha =\xi ^2 _1, \label{eq:eq6}
\\
  && { 1\over \epsilon } (d ^{2, -1} )^* \alpha + (d^{1, 0})^* \beta =
\xi _2 ^1,   \label{eq:eq7}  \\
  && (d^{1, 0})^*\alpha  + \epsilon \;  (d^{0, 1 })^* \beta =\xi ^2 _2,
\label{eq:eq8}  \end{eqnarray}
where $\xi _i = \xi _i^1 + \xi _i ^2, \xi _i^j \in \Omega ^{* , j-1}$.

(1) $l=1$. Let  $ \omega _\epsilon =\alpha  _\epsilon + \beta _\epsilon
 \in E^1_k$.
We will study the limit of $\omega _\epsilon $ as $\epsilon \to 0$.

Suppose $k\leq m-1$. By Theorem \ref{le:apr}, we see that $\Vert
\omega _\epsilon
\Vert _{H^1_{c}} $,
$\Vert \epsilon ^{-1} (d^{2, -1})^* \alpha _\epsilon \Vert _{L^2} $ and
$\Vert \epsilon^{-1}  d^{2,-1} \beta _\epsilon \Vert _{L^2}
$ are uniformly bounded. As in the proof
of Theorem 2.1,  this implies that after passing to a subsequence,
$\omega \to \omega _0 $ weakly in $H^1_c$,
where $ \omega _0 \in \Omega ^{k, 0} \bigcap Ker( (d^{2, -1}) ^* ) ={\cal R}^k
$.
So $\bar E ^1_k \subset {\cal R}^k$ for $k\leq m-1$.

We will prove that this inclusion relation also holds for $k=m$. If $k=m$,
 it follows from Theorem \ref{le:apr} that $\Vert
\beta \Vert _{H^1_c}$, $\Vert  \epsilon ^{-1}
\;  d^{2, -1} \beta_\epsilon  \Vert _{L^2}$
are uniformly  bounded. So $\beta _\epsilon
 \to 0 $ weakly in $H_c^1$
by Lemma \ref{le:d21}.

Also, since $\Vert \epsilon   (d^{1, 0}) ^* \beta _\epsilon
\Vert _{L^2} $ is  uniformly bounded
( Theorem \ref{le:apr}),
by  Lemma  \ref{le:cor}, $\epsilon \; (d^{1, 0})^* \beta
_\epsilon \to 0$ weakly
in $L^2$.
It follows from eq. (\ref{eq:eq7}) that the weak limit $\omega _0$ of
$\omega _\epsilon $ in $L^2$ satisfies
$$ ( d^{2, -1} )^* \omega _0 =0.$$
So $\bar E^1_m \subset {\cal R}^m$.

Conversely,
we will prove ${\cal R}^k \subset \bar E ^1 _k$ for  $k\leq m$.
 This follows from  the fact
that if $\alpha _0 \in \Omega ^{k, 0} / Im( d^{2, -1})$, $\Vert \alpha
_0\Vert _{L^2}=1$,  then
obviously $\alpha _0 \in E_k^1$ and hence $ \alpha _0 \in \bar E_k ^1$.

(ii) $l=2$. The proof for the case $k \ne m$ is similar to that
of Theorem 2.1 and will be omitted here.

Consider the  case $k=m$.  By Theorem \ref{le:apr}, Theorem \ref{le:apr2},
$\Vert \alpha  _\epsilon
 \Vert _{H^2_c}$, $\Vert \beta  _\epsilon \Vert _{H^1_c}$  and
$\Vert \epsilon^{-1} \; d^{2, -1} \beta  _\epsilon  \Vert _{L^2} $ are
uniformly
 bounded. Let $\alpha _0,
\beta _0$ be the weaks limit of $\alpha _\epsilon, \beta _\epsilon
 $ in $H^2_c$, $H^1_c$ respectively.
First note that, as in the case $l=1$, the weak limit $\omega_0 =\alpha
_0+\beta _0 \in {\cal R}^m$.

Moreover,  by Lemma \ref{le:cor},  $\epsilon
(d^{0, 1} ) ^* \beta  _\epsilon
 \to 0$ weakly in $L^2$. So it follows from eq. (\ref{eq:eq8})
that $(d^{1, 0})^* \alpha _0=0$.
Thus $\bar E^2_m
\subset {\cal R}^m \bigcap Ker((d^{1, 0})^*)$.

Conversely, if $\alpha _0 \in
{\cal R}^m \bigcap Ker((d^{1, 0})^*)$,
then obviously $\alpha _0 \in \bar E^2_m$.
So ${\cal R}^m \bigcap Ker($ $(d^{1, 0})^*) \subset
\bar E^2_m$.
This proves the case $l=2$.

(iii) $l=3$.  The proof in this case is similar to that of Theorem 2.1 and will
be omitted.

\end{document}